\documentclass[10pt,a4paper,reqno]{amsart} 
\usepackage{epsfig,color}
\usepackage{graphicx}
\usepackage[english]{babel} 
\usepackage{amsmath}
\usepackage{amsthm,amssymb,latexsym} 
\usepackage{t1enc}
\usepackage{fancyheadings} 
\usepackage{epic} 
\usepackage{amssymb,amsfonts,amsmath,stmaryrd}
\usepackage{eepic}

\newtheorem{prop}{Proposition} 
\newtheorem{lemma}[prop]{Lemma}
\newtheorem{corollary}[prop]{Corollary} 
\newtheorem{theorem}[prop]{Theorem}

\theoremstyle{definition}

\newcommand{\beq}{\begin{equation}}
\newcommand{\eeq}{\end{equation}}

\newcommand{\N}{\mathbb{N}} 
\newcommand{\Z}{\mathbb{Z}} 

\newcommand{\R}{\ensuremath{\mathcal{R}}} 
\newcommand{\T}{\ensuremath{\mathcal{T}}} 
\newcommand{\F}{\ensuremath{\mathcal{F}}} 
 
\newcommand{\B}{\ensuremath{\mathcal B}} 
 
\newcommand{\V}{\ensuremath{\mathcal{V}}}

\def\newop#1{\expandafter\def\csname #1\endcsname{\mathop{\rm#1}\nolimits}}
\newop{inv}

\newcommand{\Sn}{\ensuremath{\mathcal S}} 

\newcommand{\p}{permutation}
\newcommand{\ps}{permutations}
\newcommand{\gf}{generating function}
\newcommand{\gfs}{generating functions}
\newcommand{\sig}{\sigma}
\newcommand{\al}{\alpha}
\newcommand{\be}{\beta}

\newcommand{\fl}{forest-like}
\newcommand{\tl}{tree-like}

\def\emm#1,{{\em #1}}

\catcode`\@=11
\def\section{\@startsection{section}{1}%
 \z@{.7\linespacing\@plus\linespacing}{.5\linespacing}%
 {\normalfont\bfseries\scshape\centering}}
\def\subsection{\@startsection{subsection}{2}%
  \z@{.5\linespacing\@plus\linespacing}{.5\linespacing}%
  {\normalfont\bfseries\scshape}}

\def\subsubsection{\@startsection{subsubsection}{3}%
  \z@{.5\linespacing\@plus.7\linespacing}{-.5em}%
  {\normalfont\itshape}}
\catcode`\@=12
\title[Forest-like permutations]{Forest-like permutations}

\author{Mireille Bousquet-M\'elou}
\address{CNRS, LaBRI, Universit\'e Bordeaux 1,
351 cours de la Lib\'eration,
  33405 Talence Cedex, France}
\email{mireille.bousquet@labri.fr}

\author{Steven Butler}
\address{Deparment of Mathematics, University of California, San Diego, La Jolla, CA 92093-0112, USA}
\email{sbutler@math.ucsd.edu}

\date{\today}

\begin{document}

\begin{abstract}
Given a permutation $\pi\in \Sn_n$, construct a graph $G_\pi$ on
the vertex set $\{1,2,\ldots , n\}$ by joining 
$i$ to $j$ 
if (i) $i<j$ and
$\pi(i)<\pi(j)$ and (ii) there is no $k$ such that $i<k<j$ and
$\pi(i)<\pi(k)<\pi(j)$.  
We say that $\pi$ is forest-like if $G_\pi$ is a forest. We
first characterize forest-like \ps \ in terms of pattern avoidance, and
then by a certain linear map 
being onto. Thanks to recent results  of  Woo and  Yong, this  shows
that forest-like permutations characterize Schubert varieties which are
locally factorial. Thus \fl\ \ps\ generalize smooth permutations
(corresponding to smooth Schubert varieties). 

We compute the generating function of \fl \ permutations. As in the
smooth case, it turns out to be algebraic. We then adapt our method to
 count \ps\ for
which  $G_\pi$ is a tree, or a path, and recover the known \gf\ of
smooth \ps. 
\end{abstract}

\maketitle

\section{Introduction}
\begin{figure}[ht]
\begin{center}
\input{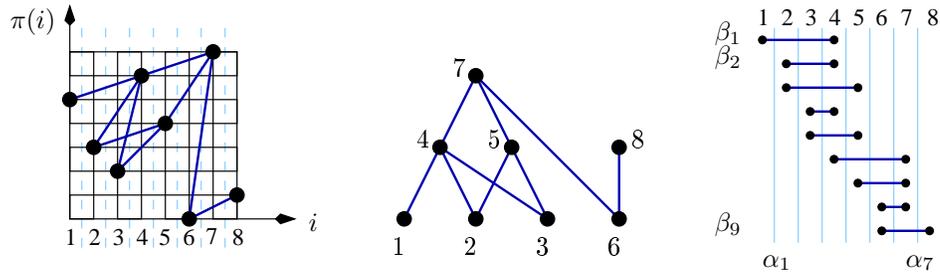}
\end{center}
\caption{The permutation $\pi= 6\ 4\ 3\ 7\ 5\ 1\ 8\ 2$, the associated
graph $G_\pi$ and the corresponding collection of bars.}
\label{fig:definitions}
\end{figure}

Take a permutation $\pi= \pi(1)\pi(2)\cdots\pi(n)$ in the symmetric
group $\Sn_n$. Let $G_\pi$ be the graph on the vertex set $\{1,2,
\ldots , n\}$ with an edge joining $i$ to $j$ 
if and only if  (i) $i<j$ and
$\pi(i)<\pi(j)$ and (ii) there is no $i<k<j$ with
$\pi(i)<\pi(k)<\pi(j)$.  An example is shown in
Figure~\ref{fig:definitions}. We say that $\pi$ is \emm \fl, if $G_\pi$ is a
forest (i.e., has no cycle).  Note that the edges of $G_\pi$
correspond to 
the edges of the \emm Hasse diagram, of the sub-poset of $\N^2$ consisting
of the points $(i,\pi(i))$ (Figure~\ref{fig:definitions}, left). This
(sub-)poset is known to play a 
crucial role in the Robinson-Schensted correspondence~\cite{greene}.

Consider also the following construction, borrowed from
\cite{WooYong1}. Label $n$ columns by $1,2, \ldots , n$, 
and place $n-1$ vertical dividers between the
columns.  
Draw a horizontal bar between column $i$ and
column $j$ if and only 
there is an edge joining $i$ and $j$ in $G_\pi$. These bars are simply
the horizontal projections of the edges of the Hasse diagram.

We use this construction to define a linear map from
$\mathbb{Z}^{n-1}$ to $\mathbb{Z}^{e(\pi)}$, where $e(\pi)$ is the
number of horizontal bars in the diagram (also 
the number of edges in $G_\pi$).  Choose a linear order on the bars,
and associate
variables $\alpha_i$ with the vertical dividers and $\beta_k$ with the
horizontal bars.  If the $k$th
horizontal bar starts in column $i$ and 
goes to column $j$ then we set 
\begin{equation}\label{eqn:linear}
\beta_k=\sum_{\ell=i}^{j-1}\alpha_\ell.
\end{equation}
 The map $L_\pi$ sends $(\al_1, \ldots, \al_{n-1})$ to $(\be_1, \ldots ,
 \be_{e(\pi)})$.  In the example above we have
$\beta_1=\alpha_1+\alpha_2+\alpha_3$, $\beta_2=\alpha_2+\alpha_3,
\ldots,\beta_9=\alpha_6+\alpha_7$. 

\medskip
Our first  result 
describes \fl\ \ps \ in terms  of the map $L_\pi$, and gives a
characterization of these \ps\ in terms of pattern avoidance. 
\begin{theorem}\label{thm:main}
For $\pi\in \Sn_n$ the following are equivalent:
\begin{itemize}
\item[(1)] the graph $G_\pi$ is a forest;
\item[(2)] the linear map $L_\pi: \Z^{n-1} \rightarrow \Z^{e(\pi)}$ is onto;
\item[(3)] the \p\ $\pi$ avoids the patterns  $1324$ and $21\bar 3 54$.
\end{itemize}
\end{theorem}

We need to clarify the third point. 
A permutation $\pi$ \emm avoids the pattern, $1324$ if one cannot find
 indices
$p<q<r<s$ such that $\pi(p)<\pi(r)<\pi(q)<\pi(s)$. 
Similarly, $\pi$ avoids   the pattern $21\bar 354$ if every occurrence
of the pattern 2154 is a subsequence of an occurrence of 21354. That is
to say, for all indices
$p<q<r<s$ such that $\pi(q)<\pi(p)<\pi(s)<\pi(r)$, there exists a
 $t$ such that $q<t<r$
and $\pi(p)<\pi(t)<\pi(s)$.  
The notation was introduced by J.~West in his
thesis~\cite{west-these}, 
and appears, for
instance, in~\cite{dulucq}.
 There are several equivalent ways to describe the latter avoidance condition.
In particular, it is easy to see that, in the terminology introduced by
Woo and Yong~\cite{WooYong1}, avoiding $21\bar 3 54$ is equivalent to
avoiding $2143$ \emm with Bruhat condition, $(1 \leftrightarrow 4)$.
However, the first description is more symmetric, more clearly showing that $\pi$
avoids $21\bar 3 54$ if and only if $\pi^{-1}$ does. 

Given that a linear map $\mathbb{Z}^{n-1} \to \mathbb{Z}^{e}$ is
bijective if and only if it is onto and $e=n-1$, we obtain the
following result.
\begin{corollary}\label{cor:main}
The map $L_\pi$ is a bijection if and only if $G_\pi$ is a tree. In
this case we say that $\pi$ is \emm tree-like.,
\end{corollary}

Our second result is the enumeration of \fl\ \ps.
We will show that their generating function is
\begin{eqnarray}
F(x)&=&{(1-x)(1-4x-2x^2)-(1-5x)\sqrt{1-4x}\over 2(1-5x+2x^2-x^3)}.
\label{Fsol}
\end{eqnarray}
We also enumerate several
natural subclasses of \fl\ \ps, such as tree-like \ps.

\smallskip
The original motivation for studying \tl\  permutations came from a
question of Woo and Yong related to Schubert varieties.  It is
known that Schubert varieties can be indexed by permutations
\cite{Fulton}, 
and various properties of Schubert varieties have been
translated into properties of permutations.  One famous example
is that a variety is smooth if and only if the associated permutation avoids the 
patterns $1324$ and $2143$~\cite{Lakshmibai}. 
A weakening of smoothness is the locally factorial property, 
an algebra-geometric
condition which states that all local rings are unique factorization
domains.  
Woo and Yong 
 established a condition for being locally
factorial which is equivalent to  $L_{\pi}$ being
onto~\cite[Prop.~2]{WooYong1}. 
They  conjectured that  this holds  if and
only if $\pi$  is $1324$ and $21\bar 3 54$
avoiding~\cite{WooYong2}. Theorem~\ref{thm:main} settles this conjecture.

We note that every \emm smooth, \p \ ($1324$ and $2143$ avoiding) is
forest-like. Smooth \ps\ have been counted before~\cite{haiman}, and
their \gf\ is:
$$
S(x)= x {\frac {  1-5\,x+4\,{x}^{2}+x\sqrt {1-4\,x}  }{1-6\,x+8
\,{x}^{2}-4\,{x}^{3}}}.
$$
As   Reference~\cite{haiman} is not easily available, we
will show how to 
adapt our proof of~\eqref{Fsol} to enumerate smooth permutations.
The series $S(x)$ occurs in several other enumeration
problems~\cite{bona-smooth}. 

\medskip
\noindent
{\bf Remark.}
%
%
Results of Cortez \cite{Cortez}, and independently Manivel
 \cite{Manivel}, show that $1324$ and $21\bar 3 54$ avoidance is 
 necessary and sufficient to characterize which
 Schubert varieties are {\em generically} locally factorial.  Here generic has 
the following sense: the variety is smooth at almost all points 
but has a closed subset $Y_{\pi}$ where it
 is not smooth, and in that closed subset it is factorial at {\em
   almost} all points.   

\medskip

We will proceed as follows. In Section~\ref{sec:characterization}
we prove Theorem~\ref{thm:main}. The proof involves a fourth
condition, equivalent to those of Theorem~\ref{thm:main}, which uses a
certain sorting procedure on the bars.
In Section~\ref{sec:generating} we count \fl\
permutations and  several of their natural subclasses, such as \tl \ \ps\
and smooth \ps.  We conclude in Section~\ref{sec:discussion} by
describing several simple bijections related to some of our enumerative
results,  and state some open problems. 

\section{Characterization of \fl\ \ps}\label{sec:characterization}
The aim of this section is to prove Theorem~\ref{thm:main}. We begin
with proving that $(1)\Rightarrow (3)$ and  $(2)\Rightarrow (3)$ by
proving the contrapositive: if $\pi$ contains  $1324$ or $21\bar 3
54$, then $G_\pi$ contains a cycle and $L_\pi$ is not onto.

\subsection{Permutations containing $1324$ or $21\bar 3 54$}
\label{sec:forbidden}
We first look at the
structure found in the diagrams of permutations containing $1324$ and
$21\bar 3 54$. 
We begin with a very simple lemma which
follows from the definition of the diagram of bars (alternatively,
from the definition of the Hasse diagram of a poset).

\begin{lemma}\label{lem:pathobars}
Let $\pi \in \Sn_n$. If $p<q$ and $\pi(p)<\pi(q)$ then there is a sequence
$p=p_0<p_1<\cdots<p_k=q$ such that $\pi(p_i)<\pi(p_{i+1})$ and in the
diagram for $\pi$ there are horizontal bars from column $p_i$ to
column $p_{i+1}$ for each $i=0,\ldots,k-1$. 
\end{lemma}

\begin{lemma}\label{lem:pattern} Given a permutation $\pi$, 
\begin{itemize}
\item[(a)]   if $\pi$ contains the pattern $1324$ then there are indices
  $p<q<r<s$ such that $\pi(p)<\pi(r)<\pi(q)<\pi(s)$ and in the diagram
  for $\pi$ there are horizontal bars from $p$ to $r$ and from $q$ to
  $s$. 
\item[(b)] if $\pi$ contains the pattern $21\bar 3 54$ then there are
  indices $p<q<r<s$ such that $\pi(q)<\pi(p)<\pi(s)<\pi(r)$ and in the
  diagram for $\pi$ there are horizontal bars from $p$ to $s$ and from
  $q$ to $r$. 
\end{itemize}
\end{lemma}
\begin{proof}
The general  idea is the following.
If we have an occurrence of the pattern that does not
satisfy  the requisite bar conditions, then
we find a  tighter occurrence that satisfies them.

For instance, start from an occurrence of  the pattern $1324$, that
is, from a sequence $p<q<r<s$ such that
$\pi(p)<\pi(r)<\pi(q)<\pi(s)$. Define $p':=\max\{i <q :
\pi(i)<\pi(r)\}$ and $r':=\min\{j>q : \pi(r) \ge \pi(j)>\pi(p')\}$. Then
$p\le p'<q<r'\le r<s$, the sequence $p', q, r', s$ corresponds to another
occurrence of $1324$, and there is a bar between columns $p'$ and $r'$. 

The rest of the lemma is proved by similar arguments.
\end{proof}

\noindent{\bf Remark.} Point (b) in the above lemma shows that
  $21\bar 3 54$ avoidance can also be described graphically as follows. Take a
  permutation $\pi$ and plot it as on the left of
  Figure~\ref{fig:definitions}. Represent by straight lines
  the edges of the Hasse diagram of the poset $\{(i, \pi(i))\}$. We thus
  obtain the \emm natural embedding, of $G_\pi$. Then
  $\pi$ avoids $21\bar 3 54$ if and only if  \emm this embedding, of 
$G_\pi$ is planar (no edges cross).  This does not mean that avoiding
  $21\bar 354$ is necessary for $G_\pi$ to be planar: for
  instance, the \p\ $\pi=2143$ contains $21\bar354$ but $G_\pi$ is
  planar (though its natural embedding is not). 

\medskip

\begin{figure}[ht]
\begin{center}
\input{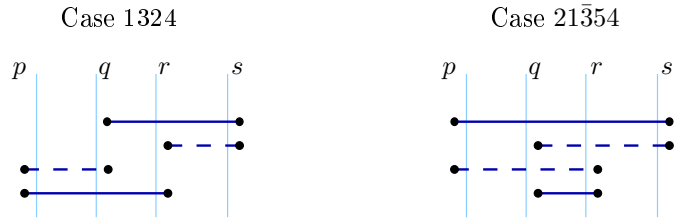}
\end{center}
\caption{Patterns in the bar diagrams of \ps\ containing 1324 or
  $21\bar 3 54$.} 
\label{fig:forbidden}
\end{figure}

Lemma~\ref{lem:pattern} is illustrated in Figure~\ref{fig:forbidden},
where the solid lines 
indicate a single bar and the dashed lines indicate a sequence of
bars (coming from Lemma~\ref{lem:pathobars}). 

We now show that the occurrence of either of the two ``forbidden'' patterns
  implies the existence of cycles in $G_\pi$, and prevents $L_\pi$
  from being onto.
First, from Figure~\ref{fig:forbidden} we can read off cycles in
$G_\pi$.  For example in 
the $1324$ case we have a cycle that starts at $p$, goes to
$r$ then by a sequence of edges goes to $s$ then to $q$
and finally by another sequence of edges we return to $p$.  This is a
true cycle, as it contains the edge joining $p$ to $r$ only once.
Similarly, in the $21\bar 3 54$ case, there is a true cycle visiting $p, s,
q, r$ in this order. 
Secondly, we also see that there  are nontrivial linear
dependencies among the $\beta_j$.  In the $21\bar 3 54$ case the sum
of the solid bars equals the sum of the dashed bars, and a similar event
happens in the $1324$ case.  This prevents $G_\pi$ from being onto.

So if the permutation contains $1324$ or
$21\bar 3 54   $ then $G_\pi$ has cycles and $L_\pi$ is not onto.
Taking the contrapositive 
gives (1)$\Rightarrow$(3) and (2)$\Rightarrow$(3) in Theorem~\ref{thm:main}.

\subsection{Sorting the horizontal bars}
\label{sec:sorting}
In this subsection, we define a new condition (2') that clearly
implies the  surjectiveness condition (2). We then prove that (2') is 
implied by the pattern avoidance condition (3), and finally that (2')
implies the acyclicity condition (1). Combined with
Section~\ref{sec:forbidden}, this proves that the four
conditions (1), (2), (2') and (3) are equivalent, and establishes 
Theorem~\ref{thm:main}.
The structure of the proof is schematized below.

\smallskip
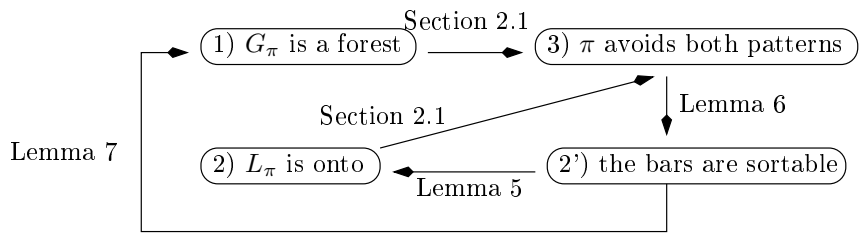
\begin{center}
  \input{schema.pstex_t}
\end{center}

In the construction of the diagram for a permutation we placed no
condition on the ordering of the horizontal bars from top to bottom.
We now describe a way to attempt to sort them.
  Create a second diagram with the
same columns but no horizontal bars.  We now look for bars to move to
the second diagram by scanning the vertical dividers 
from left to right,
 looking for any
divider which is intersected by exactly one horizontal bar.  As soon
as we find such an intersection we move the corresponding horizontal bar to
the second diagram and put it above
 any previously moved bar.  We
then repeat this scanning process, starting again from the leftmost
divider, until no divider intersects exactly 
one horizontal bar.  If at this stage all the horizontal bars
are moved over, we say that
\begin{quote}
 (2') \emm the bars are fully sortable.,  
\end{quote}

\begin{figure}[ht]
			\begin{center}
\input{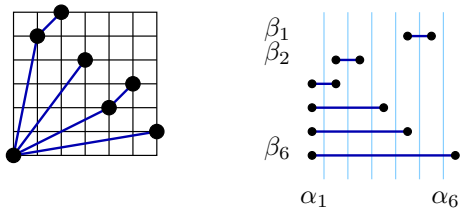}
\end{center}
\caption{The permutation
  $\pi= 1\ 6\ 7\ 5\ 3\ 4\ 2$ and the associated sorted diagram of bars.}
\label{fig:sorting}
\end{figure}

\noindent By construction, this can only happen when the
number of edges satisfies $e(\pi)\le n-1$. 
An example of a fully sorted diagram is shown
in Figure~\ref{fig:sorting}. 

Assume the bars are fully sortable. In terms of the
equations~\eqref{eqn:linear}, this means that \emm at least one new
variable, $\al_i$ \emm occurs in each equation,.  More precisely, if
$\V_k$ denotes the set of variables $\al_i$ occurring in $\be_1,
\ldots , \be_k$, then $\V_k \subsetneq \V_{k+1}$.  Hence, given $\be \in
\Z^{e(\pi)}$, the system~\eqref{eqn:linear} can be solved for $\al$
by backward 
substitution  from the top equation to the bottom
equation. Consequently, we have the following.

\begin{lemma}\label{lem:onto}
If the bars are fully sortable then $L_\pi$ is onto.
\end{lemma}

In other words,  (2') implies (2). 
We shall see below that the converse is also true. 
This will be a consequence of Theorem~\ref{thm:main} and 
the following lemma, which proves (the contrapositive of) (3)$\Rightarrow$(2').

\begin{lemma}\label{lem:failsort}
If  the bars are not fully sortable then $\pi$ contains $1324$ or
$21\bar 3 54$. 
\end{lemma}
\begin{proof}
If we stopped before all the bars have been moved over then it must be
the case that for what remains all the vertical dividers intersect
either zero, or two or more horizontal bars.  We will work with these
remaining (i.e., unmoved) horizontal bars. 

Suppose that column $a$ is the leftmost column which has the start of
a bar, then as noted above it must be the start of at least two bars
(otherwise we would have moved the bar over).  Let $c$ denote the
column where the {\em longest} horizontal bar starting in column $a$
ends.  Let $b$ be the 
 rightmost column satisfying $a<b<c$ and
$\pi(c)<\pi(b)$
 (such a $b$ exists because the end of a second bar
 that starts in $a$ satisfies both conditions). 

We now consider cases on how to cover the vertical divider to the
right of column $b$ with a second horizontal bar.   

{\em Case $(1)$.} There is a horizontal bar that begins at $b$. This bar
 ends at  some position $d$, 
which, by the choice of $b$, satisfies $d>c$. 
 In this case we have that $a<b<c<d$
while $\pi(a)<\pi(c)<\pi(b)<\pi(d)$ and so $\pi$ contains the pattern
$1324$. 

{\em Case $(2)$.} There is a horizontal bar that begins at column $d$
where $d<b$ and crosses to some column $e$ where $e>b$.  
By the choice of $b$, we have $a<d$. 
Since $d$
lies between $a$ and $c$ we must have that $\pi(d)<\pi(a)$ or
$\pi(d)>\pi(c)$
(if $\pi(d)$ were in the interval $[\pi(a), \pi(c)]$, there would not
be a bar from $a$ to $c$).
  So we consider subcases. 

{\em Case $(2a)$.} If $\pi(d)<\pi(a)$ then we have that $a<d<b<c$ and
$\pi(d)<\pi(a)<\pi(c)<\pi(b)$ and since there is a horizontal bar from
$a$ to $c$, $\pi$ contains the pattern $21\bar 3 54$.  (Note this includes
the possibility that $c=e$.) 

{\em Case $(2bi)$.} Suppose that not only $\pi(d)>\pi(c)$, but also 
$\pi(d)>\pi(b)$.  Then we note that we
have $a<d<b<e$ and $\pi(a)<\pi(b)<\pi(d)<\pi(e)$ and so $\pi$ contains the
pattern $1324$. 

{\em Case $(2bii)$.} Suppose finally that $\pi(c)<\pi(d)<\pi(b)$. By the
choice of $b$, we must have $e>c$.
Then we note that we have $a<d<c<e$ and $\pi(a)<\pi(c)<\pi(d)<\pi(e)$
and so $\pi$ contains the pattern $1324$. 
\end{proof}

Our final lemma  proves  that (2')$\Rightarrow$(1). 
\begin{lemma}\label{lem:acyclic}
  If the bars are fully sortable, then $G_\pi$ is a forest.
\end{lemma}
\begin{proof}
 Suppose on the contrary that $G_\pi$ contains a cycle and we can
 fully sort the bars.  Now consider the set $\B$ of bars that
correspond to the edges of a cycle in $G_\pi$. At some stage in the sorting
procedure, a first bar $b$ of $\B$ is moved over. At this stage,
it is the only bar that crosses some vertical divider, say, the $i$th
one. In particular, all the other bars involved in the cycle
lie entirely  to the right or entirely to the left of the $i$th divider.
In terms of $G_\pi$, this means that removing the edge corresponding
to $b$ has \emm disconnected,  the cycle. This is of course
impossible, so $G_\pi$ cannot contain a cycle.
 \end{proof}

\section{Generating functions for forest-like permutations}
\label{sec:generating}
We now want to prove the enumerative result~\eqref{Fsol}. 
At the heart of this result is a recursive description of \fl\ \ps,
given in Proposition~\ref{prop:decompose}. This decomposition is then
translated into a functional equation defining the \gf\ of \fl\ \ps\
(Proposition~\ref{prop:bivariate}), which we solve using the 
\emm kernel method,. 

The same decomposition can be recycled to count various subclasses of
\fl\ \ps.
We will thus also obtain the \gfs\ of
\begin{enumerate}
\item \tl \ \ps , 
\item  \emm rooted, \tl\ \ps\ (the term \emm rooted,
meaning that $\pi(1)=1$),
\item \emm path-like, \ps\  ($G_\pi$ is a path),
\item  \emm smooth, \ps\ ($\pi$ avoids 1324 and 2143). 
\end{enumerate}
Note that every \fl \ \p\ satisfying $\pi(1)=1$ is actually \tl\
(every vertex of $G_\pi$ is connected to the vertex $1$), and thus is
a rooted \tl\ \p.
 Note also the following inclusions:
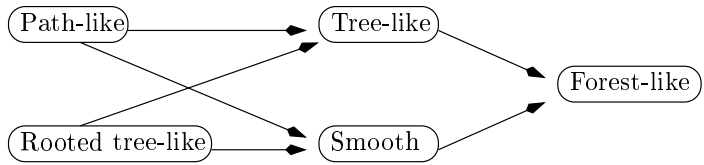
\begin{center}
  \input{inclusions.pstex_t}
\end{center}

For $n\ge 1$, we denote by $f_n$ (resp.~$t_n, r_n, p_n, s_n$) the  number of
permutations $\pi\in \Sn_n$ of the above five types. We introduce the
corresponding \gfs \ $F(x)$ (resp.~$T(x), R(x), P(x), S(x)$). In
particular,
$$
F(x)=\sum_{n\geq
  1}f_nx^n= x+2x^2+6x^3+22x^4+89x^5+379x^6+1661x^7+\cdots.
$$

Our enumerative results are summarized in the following theorem.
\begin{theorem}\label{thm:enumerate}
The five generating functions defined above are given by:
\begin{eqnarray*}
F(x)&=&{(1-x)(1-4x-2x^2)-(1-5x)\sqrt{1-4x}\over
  2(1-5x+2x^2-x^3)},\\
T(x)&=& 
{1-3x-6x^2-(1-5x)\sqrt{1-4x}\over
  2(2-9x)},\\
R(x)&=&{1-\sqrt{1-4x}\over2},\\
P(x)&=& x\,\frac {1-2x+2x^2}{(1-x)(1-2x)},\\
S(x)&= & x {\frac {  1-5\,x+4\,{x}^{2}+x\sqrt {1-4\,x}  }{1-6\,x+8
\,{x}^{2}-4\,{x}^{3}}}.
\end{eqnarray*}
\end{theorem}

From these generating functions it can be shown that there exists
positive constants $\kappa$ such that
$$
  f_n\sim\kappa_f(4.61\ldots)^n, \qquad
t_n\sim\kappa_t(4.5)^n,\qquad
r_n={1\over n}{2n-2\choose n-1}\sim\kappa_r 4^{n-1}n^{-3/2},
$$
$$
 p_n=2^{n-1}-1\  \hbox{ for } n\ge 2,\qquad \qquad 
s_n \sim \kappa_s (4.38\ldots)^n,
$$
where the growth
constants occurring in the asymptotics of $f_n$ and $s_n$ are
respectively  the real roots of the polynomial $t^3-5t^2+2t-1$ and
$t^3-6t^2+8t-4$. 

We note that $r_n$ is the $(n{-}1)$st Catalan number and has numerous
combinatorial interpretations~\cite[Chap.~6]{stanley-vol2}.  
We give in Section~\ref{sec:plane} a bijective proof of this result,
as well as another bijection explaining why the numbers $p_n$  are so
simple. The terms $t_n$ 
have also occurred before and enumerate the number of \emm stacked
directed animals on a triangular lattice, \cite{mbm-rechni}.  No
direct bijection between stacked directed animals and tree-like
permutations is currently known.  

\medskip
The form of our decomposition of \fl\ \ps\ will force us to
 take into account an additional statistic, namely the number
of $rl$-minima for \fl\ or \tl\ \ps, and the length of the final
ascent in smooth \ps. This is why we actually obtain bivariate \gfs\
that refine the above theorem
(see~\eqref{eq:T},~\eqref{eq:F},~\eqref{eq:S}). Other statistics,
like the number of descents, could also 
be carried through our calculations.

\subsection{Decomposing forest-like permutations}
If $\pi\in \Sn_n$ we say that $\pi$ has length $n$, and write
$|\pi|=n$.  We say that $\pi(i)$ is an $rl${\em-minimum\/}
(right-to-left-minimum) if for all $j>i$, we have $\pi(j)>\pi(i)$.
We denote by $m(\pi)$ the number of $rl$-minima of $\pi$. Finally,
$\pi$ is \emm increasing, if $\pi=12\cdots n$.

Let $\pi \in \Sn_n$ be \fl. We decompose $\pi$ by considering which
element maps to $1$. 
So suppose that $i=\pi^{-1}(1)$ then there are two cases:  

\noindent $\bullet$ \emm First case,: $i=n=|\pi|$. Then the permutation $\tau\in
\Sn_{n-1}$ defined by $\tau(i)=\pi(i)-1$ is
forest-like.  Conversely, starting with a forest-like permutation
$\tau\in \Sn_{n-1}$ we can construct a \fl\ permutation $\pi\in \Sn_n$ by
letting $\pi(i)=\tau(i)+1$ for $1\leq i\leq n-1$ and $\pi(n)=1$.  Note
that  $\pi$  is tree-like if and only if $n=1$. 

\begin{figure}[htb]
\begin{center}
  \input{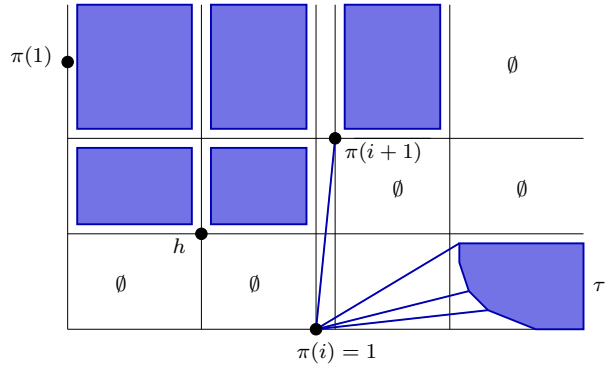}
\end{center}
\caption{The structure of a forest-like \p . The shaded areas show which
regions of the embedding in $\mathbf{N}^2$ may contain points.}
%
%
\label{fig:structure}
\end{figure}

\noindent $\bullet$ \emm Second case,: $i=\pi^{-1}(1)<n$.
We now focus on this case, illustrated in Figure~\ref{fig:structure}.  Let  
\beq\label{h-def}
h=\min\big(\{\pi(i+1)\}\cup\{\pi(j):j<i\big\}).
\eeq
So $h$ is the smaller of the lowest value of $\pi$ to the left of $i$
or the value of $\pi$ at $i+1$. 

First note that for all $j\geq i+1$ we have $\pi(j)\leq h$ or $\pi(j)\geq
\pi(i+1)$.  If not, then for some $j$ we have
$\pi(i)<h<\pi(j)<\pi(i+1)$ and $\pi^{-1}(h)<i<i+1<j$, so the
permutation contains the pattern $21\bar 354$, and  cannot be 
forest-like.  Further, if $j,k\geq i+1$ with $\pi(j)\geq \pi(i+1)$ and
$\pi(k)<h$ then $j<k$.  If not, then $i<i+1<k<j$ and
$\pi(i)<\pi(k)<\pi(i+1)<\pi(j)$, so the permutation contains the pattern
$1324$, and  cannot be  forest-like. 

The latter property implies that the last $h-2$ terms of
$\pi=\pi(1)\pi(2) \cdots \pi(n)$ are
$2,3,\ldots,h-1$ in some order.  Let $\tau$ be the 
permutation obtained from $\pi$ by retaining only its $h-1$ smallest
entries, i.e., 
\[
\tau=1\,\pi(n-h+3)\,\pi(n-h+4)\cdots\pi(n).
\]
Then $\tau$ is rooted and \tl. Similarly, let $\sigma$ be the
 permutation obtained by deleting these $h-1$ smallest
entries and subtracting $h-1$ from the remaining entries:
\[
\sig=\big(\pi(1)-h+1\big)\cdots\big(\pi(i-1)-h+1\big)\big(\pi(i+1)-h+1\big)\cdots\big(\pi(n-h+2)-h+1\big).
\]
Then $\sigma$ is \fl. Moreover, 
$\sigma(i)$ is  an $rl$-minimum of $\sigma$. If $\sigma(i)$ is the
$k$th $rl$-minimum of $\sigma$ (read \emm from  right to left,\,), define
$\Phi(\pi)=(\tau,\sigma,k)$. Observe that $k=m(\sigma)$ if
$\sigma(i)=1$ (that is to say, $h=\pi(i+1)$), and $k=1$ if  $i=|\sigma|$.

Conversely, starting from a 3-tuple $(\tau,\sigma,k)$ such that $\tau$
is \tl, $\sigma$ is \fl\ and $k\le m(\sigma)$, we can construct a
(unique) \fl\ \p\  $\pi$ satisfying $\Phi(\pi)=(\tau,\sigma,k)$.
If $|\tau|=h-1$ and the $k$th $rl$-minimum of $\sigma$ is $\sigma(i)$,
this is done by adding $h-1$ to the entries of $\sigma$, inserting $1$ 
at position $i$ and adding the 
other entries of $\tau$ to the right of $\sigma$, in the same order as
in $\tau$.  
By looking at the number of $rl$-minima of the resulting \p\ $\pi$,
we obtain the following result.  

\begin{prop}\label{prop:decompose}
The map $\Phi$ is a bijection between forest-like permutations $\pi$
with $\pi^{-1}(1)<|\pi|$ and $3$-tuples $(\tau,\sigma,k)$ such that
$\tau$ is rooted tree-like, $\sigma$ is forest-like, and $1\leq k\leq
m(\sig)$.  Moreover, 
\beq\label{m-rec}
|\pi|=|\tau|+|\sigma|\qquad\mbox{and}\qquad
m(\pi)=\left\{\begin{array}{l@{\quad}l}  k+1
&\mbox{if }\tau= 1,\\
m(\tau)&\mbox{otherwise}.\end{array}\right.
\eeq
\end{prop}
In order to count the various sub-classes of \fl\ \ps\ we have
defined, we need the following result.
\begin{prop}\label{prop:subclasses}
Let $\pi=\Phi(\tau,\sigma,k)$ be a \fl\ \p\ such that $\pi^{-1}(1)<|\pi|$. Then
\begin{enumerate}
\item  $\pi$ is tree-like if and only if $\sigma$ is tree-like, 
\item  $\pi$ is rooted tree-like if and only if $\sigma$ is rooted tree-like
and $k=m(\sigma)$, 
\item $\pi$ is path-like if and only if $\tau$ is increasing, 
$\sigma$  is   path-like and its $k$th $rl$-minimum $\sig(i)$ is
  such that $i$ has degree $1$ in $G_\sig$,
\item $\pi$ is smooth if and only if $\sigma$ is smooth and either
  $k=m(\sigma)$ or $k\le  a(\sigma)$, where $a(\sigma)$ is the
  \emm length of the final ascent, of $\sigma$: if $|\sigma|=\ell$,
\beq\label{a-def}
a(\sigma)=\max\{i : \sigma(\ell-i+1)< \cdots <\sigma
(\ell-1)<\sigma(\ell)\}.
\eeq
 Moreover, 
\beq\label{a-rec}
a(\pi)=\left\{\begin{array}{l@{\quad}l} 
k+1 &\mbox{if }\tau= 1 \mbox{ and } k\le a(\sigma),\\
a(\sigma) &\mbox{if }\tau= 1 \mbox{ and } k=m(\sigma)> a(\sigma),\\
a(\tau)-1&\mbox{if }  \tau \not = 1 \mbox{ is increasing},\\
a(\tau)&\mbox{otherwise}
.\end{array}\right.
\eeq
\end{enumerate}
\end{prop}
\begin{proof}
The first three results simply follow from the decomposition of
  Figure~\ref{fig:structure}. The reader should look at Figures~\ref{fig:plane}
  and~\ref{fig:path-decomp} to see this decomposition specialized to the
  rooted case   and the path case, respectively. 

Now let's assume that $\pi$ is
  smooth. Since $\sigma$ and $\tau$ are obtained by deleting entries
  from $\pi$, they are smooth as well. This does not restrict the
  choice of $\tau$, since  every rooted tree-like \p\ is
  smooth. 
Conversely, when we
  construct $\Phi(\tau,\sigma,k)$ (assuming that $\sigma$ is smooth
  and $\tau$ rooted) we do not create any occurrence
  of $2143$ if we insert 1 just before the smallest entry of
  $\sigma$. This corresponds to the case $k=m(\sigma)$. 

However, if
  $k<m(\sigma)$, 
  then the value $h$ defined by~\eqref{h-def} satisfies $h<\pi(i+1)$, and
  the final permutation $\pi$ contains 2143 if, and only if, there is
  a descent in $\sigma$ somewhere to the right of $\sigma(i)$. In
  other words,  if
  $k<m(\sigma)$, then $\pi$ avoids 2143 if and only if 1 is inserted
  in the final ascent of $\sigma$, that is to say, $k\le a(\sigma)$.

A case study finally provides the value of $a(\pi)$.
\end{proof}

\subsection{Functional equations}
We now translate Propositions~\ref{prop:decompose}
and~\ref{prop:subclasses} into enumerative terms.  We 
first note that every pair $(\tau,\sigma)$ can be combined in
$m(\sigma)$ different ways. To account for this we refine our
\gfs\  by further distinguishing by the number of $rl$-minima.
So let 
$$
\mathcal{F}(u)\equiv \mathcal{F}(x,u)=\sum_{n,\ell\geq
  1}f_{n,\ell}\,x^nu^{\ell}=\sum_\ell\mathcal{F}_\ell(x)u^\ell
$$
where $f_{n,\ell}$ is the number of forest-like permutations of $\Sn_n$
having $\ell$ 
 $rl$-minima. Note that
$F(x)=\mathcal{F}(1)$. Define similarly the bivariate series
$\mathcal{T}(x,u)$, 
$\mathcal{R}(x,u) $, $\mathcal{P}(x,u) $. 
The case of smooth permutation is a bit different: here, the crucial
parameter is the length of the final ascent, defined
by~\eqref{a-def}.  We thus use a new indeterminate $v$ and define
$$
\mathcal{S}(v)\equiv \mathcal{S}(x,v)=\sum_{n,\ell\geq
  1}s_{n,\ell}\, x^nv^{\ell}=\sum_\ell\mathcal{S}_\ell(x)v^\ell
$$
where $s_{n,\ell}$ is the number of smooth permutations of $\Sn_n$
having a final ascent of length $\ell$. We define similarly the series
$\overline{\mathcal{R}}(x,v)$ that counts rooted tree-like \ps\ by the
same statistics.

\begin{prop}\label{prop:bivariate}
The (bivariate) generating functions $\mathcal{F}(u)$,
$\mathcal{T}(u)$, $\mathcal{R}(u)$ and $\mathcal{P}(u)$ satisfy:
\[
\begin{array}{r@{~=~}l@{~+~}c@{~+~}l}
\mathcal{F}(u)&xu~+~xu\mathcal{F}(1)&\displaystyle xu^2{\mathcal{F}(u)-\mathcal{F}(1)\over u-1}&(\mathcal{R}(u)-xu)\mathcal{F}'(1),\vspace{5pt}\\
\mathcal{T}(u)&xu&\displaystyle xu^2{\mathcal{T}(u)-\mathcal{T}(1)\over u-1}&(\mathcal{R}(u)-xu)\mathcal{T}'(1),\vspace{5pt}\\
\mathcal{R}(u)&xu&xu\mathcal{R}(u)&(\mathcal{R}(u)-xu)\mathcal{R}(1),\vspace{5pt}\\ 
 \mathcal{P}(1)&x& \displaystyle\frac{x^2}{(1-x)^2} & \displaystyle\frac x{1-x} \left(
 \mathcal{P}(1)-x\right) ,
\end{array}
\]
where $\mathcal{F}'(1)={\partial \mathcal{F} \over \partial u}(x,1)$ and
similarly for $\mathcal{T}'(1)$.  Moreover, 
\begin{equation}\label{link}
\mathcal{F}(u)={\mathcal{T}(u)\over1-\mathcal{T}(1)}.
\end{equation}
For the smooth case,
\begin{multline*}
  \mathcal{S}(v)=xv(1-x)+x\mathcal{S}(v)+xv(1-x)\frac{v\mathcal{S}(v)-\mathcal{S}(1)}{v-1}\\
+\left( \overline{\mathcal{R}}(v)-\frac{xv(1-x)}{1-xv}\right) \left(
(1-x)(\mathcal{S}'(1)+\mathcal{S}(1))-x\right)
\end{multline*}
where
\[
\mathcal{\overline{\mathcal{R}}}(v)= \frac{xv(1-x)}{1-xv} +x\overline{\mathcal{R}}(v) + \left(
\overline{\mathcal{R}}(v)-\frac{xv(1-x)}{1-xv}\right) \overline{\mathcal{R}}(1).
\]

\end{prop}

\begin{figure}[htb]
\begin{center}
  \input{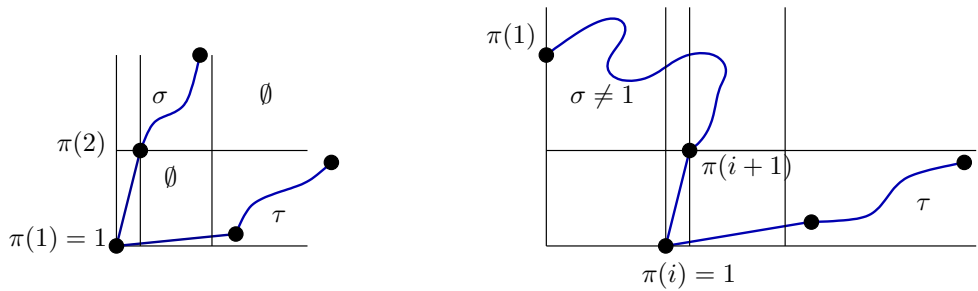}
\end{center}
\caption{The decomposition of path-like \ps.}
\label{fig:path-decomp}
\end{figure}

\begin{proof}
We first consider $\mathcal{F}(u)$.  The terms $xu+xu\mathcal{F}(1)$
count forest-like permutations with $\pi^{-1}(1)=|\pi|$, which have
only one $rl$-minimum.  For the remaining forest-like permutations we
use Proposition~\ref{prop:decompose}.  The \gf\ of \ps\ $\sigma$ such
that $\tau=1$ is:
\[
\sum_\ell\mathcal{F}_\ell(x)\sum_{k=1}^\ell
xu^{k+1}~=~xu^2\sum_\ell\mathcal{F}_\ell(x){u^\ell-1\over
  u-1}~=~xu^2{\mathcal{F}(u)-\mathcal{F}(1)\over u-1}, 
\] 
while for the \ps\ such that  $\tau\neq 1$ we obtain:
\[
\sum_\ell\mathcal{F}_\ell(x)\sum_{k=1}^\ell\big(\mathcal{R}(u)-xu\big)~=~\big(\mathcal{R}(u)-xu\big)\mathcal{F}'(1). 
\]
Combining all cases gives the result for $\mathcal{F}(u)$.  

The equation for $\mathcal{T}(u)$ is proved in a similar way  (note
that there is no counterpart to the term 
$xu\mathcal{F}(1)$ since this corresponds to forests where $1$ is an
isolated vertex). 

For rooted tree-like permutations there is no choice
in the way we merge $\tau$ and $\sigma$ and so we obtain a
significantly simpler equation (see Figure~\ref{fig:plane}). 

\medskip

The equation we have obtained for $\F(u)$ shows that the
indeterminate $u$ is needed to exploit the decomposition of
Proposition~\ref{prop:decompose}. This is not the case for path-like
\ps, and this is why we will not take into account the number of
$rl$-minima. 
If $\sig$ is path-like, the graph $G_\sig$ has exactly 2
vertices of degree 1, unless $\sig=1$. If $\sig$ is increasing, both of these end
vertices correspond to $rl$-minima. Otherwise, only the largest one
does (Figure~\ref{fig:path-decomp}). The term $x^2/(1-x)^2$ in the equation corresponds to the case
where $\sig$ is increasing and $k=|\sig|$. The term
$x/(1-x)(\mathcal{P}(1)-x)$ corresponds to the case $k<|\sig|$. 

The relationship~\eqref{link} can be explained by noting that a
forest-like permutation $\pi$ is either tree-like, or is obtained by
appending a tree-like permutation $\tau$ to the beginning of another
forest-like permutation $\sig$. More formally,  
\[
\pi=\big(\tau(1)+h\big)\big(\tau(2)+h\big)\cdots\big(\tau(k)+h\big)\sigma(1)\sigma(2)\cdots\sigma(h),
\]
where $\tau$ is tree-like and $\sigma$ is forest-like.  Note that
$m(\pi)=m(\sigma)$.  In terms of \gfs, this gives
$\mathcal{F}(u)=\mathcal{T}(u)+\mathcal{T}(1)\mathcal{F}(u)$. 

\medskip

We now proceed with the smooth case. Let us first determine the
\gf\ $\mathcal{S}_0(v)$ counting the smooth
\ps\ $\pi$ such that $a(\pi)=m(\pi)$ (that is to say, 1 belongs to the
final ascent of $\pi$). This equality certainly holds if
$\pi^{-1}(1)=|\pi|$. Otherwise, let us write
$\pi=\Phi(\tau,\sigma,k)$. By comparison of~\eqref{m-rec}
and~\eqref{a-rec}, we see that  $a(\pi)=m(\pi)$ if and only if
$\tau=1$ and $k\le a(\sig)$.  Hence
\beq\label{S0}
\mathcal{S}_0(v)= xv(1+\mathcal{S}(1))+ x\sum_\ell \mathcal{S}_\ell(x)
\sum_{k=1}^\ell v^{k+1}= 
xv + xv \,\frac{v\mathcal{S}(v)-\mathcal{S}(1)}{v-1}.
\eeq
Combining this with~\eqref{a-rec}, it follows that the smooth \ps\
$\pi=\Phi(\tau,\sigma,k)$ such that 
$\tau=1$ but $k=m(\sig)>a(\sig)$  are counted by
\beq\label{S1}
x\left( \mathcal{S}(v)-\mathcal{S}_0(v)\right) .
\eeq
In the case where $\tau\not=1$ is increasing, we obtain the series
\beq\label{S2}
\frac{x^2v}{1-xv} \left( \mathcal{S}'(1)+\mathcal{S}(1)-\mathcal{S}_0(1)\right)
\eeq
while in the case  where $\tau$ is not increasing, we find:
\beq\label{S3}
\left(\overline{\mathcal{R}}(v)-\frac{xv}{1-xv} \right) \left( \mathcal{S}'(1)+\mathcal{S}(1)-\mathcal{S}_0(1)\right).
\eeq
The series $\mathcal{S}(v)$ is the sum of~(\ref{S0}--\ref{S3}). This gives the
desired
 functional equation for $\mathcal{S}(v)$.

It remains to count rooted \tl\ \ps\ by the length of the final
ascent. We obtain an equation for $\overline{\mathcal{R}}(v)$ by specializing the above
study to the rooted case, that is to say, to the case where $\sigma$
is rooted and $k=m(\sigma)$. The counterparts of
the terms~(\ref{S0}--\ref{S3}) are respectively
\begin{multline*}
  \overline{\mathcal{R}}_0(v)= \frac{xv}{1-xv}, ~~ x\left(\overline{\mathcal{R}}(v)-\overline{\mathcal{R}}_0(v)\right), ~~
\frac{x^2v}{1-xv} 
{\mathcal{\overline{\mathcal{R}}}}(1) \  \mbox{ and }
\left(\overline{\mathcal{R}}(v)-\frac{xv}{1-xv}
\right){\mathcal{\overline{\mathcal{R}}}}(1).
\end{multline*}
The sum of these four terms is $\overline{\mathcal{R}}(v)$, and this
gives the desired equation.
\end{proof}

\subsection{Solution of the functional equations}
We are finally going to solve the equations of
Proposition~\ref{prop:bivariate} to obtain
Theorem~\ref{thm:enumerate}. Three of them do not raise any
difficulty. 
Namely, the equation defining $\mathcal{P}(1)$ is readily solved,
while the equations defining  $\mathcal{R}(u)$ 
and  $\overline{\mathcal{R}}(v)$ can be solved by first setting $u=1$
(or $v=1$) to determine the value of these series at $u=1$ (or $v=1$)
and then using these preliminary results  to compute the full series.
In particular,
\begin{equation}\label{eq:R}
\mathcal{R}(u)~=~{xu(2-u-u\sqrt{1-4x}) 
\over 2(1-u+xu^2)}~=~{xu\over 1-u\mathcal{R}(1)}.
\end{equation}

The other three equations (defining $\R$, $\T$ and $\mathcal S$) involve
divided differences of the form  
$$
\frac{A(u)-A(1)}{u-1}
$$
and cannot be solved by setting $u=1$.  Instead, we will solve them by
using the \emm kernel method,~\cite{hexacephale,bousquet-petkovsek-1}.
Consider for instance the 
equation for tree-like permutations.  This is a linear 
equation with one \emm catalytic, variable ($u$) and two additional unknown
functions ($\mathcal{T}(1)$ and $\mathcal{T}'(1)$).  However, these
two functions are not independent: by taking the limit as
$u$ goes to $1$ in the equation we find 
\beq\label{T1T1p}
\mathcal{T}(1)~=~x+\mathcal{R}(1)\mathcal{T}'(1).
\eeq
The coefficient of $\mathcal{T}(u)$ in the equation defining $\T(u)$ is
\[
1-{xu^2\over u-1}~=~{u-1-xu^2\over u-1},
\]
which vanishes for two values of $u$.  One of these values is a formal
power series in $x$, 
\[
U~\equiv~U(x)~=~{1-\sqrt{1-4x}\over 2x}.
\]
Replacing $u$ by $U$ in the functional equation gives a second linear
relation between $\mathcal{T}(1)$ and $\mathcal{T}'(1)$:
\beq\label{T1T1p2}
0=xU-\mathcal{T}(1)+\big(\mathcal{R}(U)-xU\big)\mathcal{T}'(1).
\eeq
One can now solve~\eqref{T1T1p} and~\eqref{T1T1p2} for
$\mathcal{T}(1)$ and $\mathcal{T}'(1)$, in terms of $x, U,
\mathcal{R}(1)$ and $ \mathcal{R}(U)$. Then the solution can be written as
a pair of rational functions of $U$ using:

-- the expression of   $\mathcal{R}(U)$ in terms of
$x, U$ and $\mathcal{R}(1)$ (see~\eqref{eq:R}), 

--  the fact that
$\mathcal{R}(1)=xU$, 

-- the equation $x=(U-1)/U^2$.

\noindent 
Replacing the expressions of $\mathcal{T}(1)$ and $\mathcal{T}'(1)$ in
the original functional equation  gives an
expression for $\mathcal{T}(u)$ in terms of $u$ and $U$, which can be
rewritten as
\begin{equation}\label{eq:T}
\mathcal{T}(x,u)~=~xu{(1+V)^2(1-2V)-uV(1-2V-2V^2)\over(1-2V)(1+V-uV)^2}
\end{equation}
where
\[
V~=~U-1~=~{1-2x-\sqrt{1-4x}\over 2x}.
\]

We can use similar techniques to find $\mathcal{F}(u)$. However,
it is easier to use~\eqref{link} and what we
have obtained 
for $\mathcal{T}$ to get 
\begin{equation}\label{eq:F}
\mathcal{F}(x,u)~=~uV{(1+V)^2(1-2V)-uV(1-2V-2V^2)\over(1-V-2V^2-V^3)(1+V-uV)^2}
\end{equation}
where $V$ is given above.

The solution of the equation defining $\mathcal {S}(u)$ is similar to
what we have done for $\T(u)$. One possible expression of the
bivariate series that counts smooth \ps\ by the length and the length
of the final ascent is
\beq\label{eq:S}
\mathcal{S}(x,u)=xu\, {\frac { 
  \left( 1+V \right)  \left( 1-{V}^{2}-{V}^{3}
 \right) -Vu \left( 1-V-{V}^{2}-{V}^{3}   \right) }{ \left( 1+
V-uV \right)  \left( 1-V-{V}^{2}-{V}^{3} \right)  \left( 1-xu
 \right) }}.
\eeq
Putting $u=1$ into equations \eqref{eq:T}, \eqref{eq:F} and
\eqref{eq:S} and simplifying 
then gives the results of Theorem~\ref{thm:enumerate}. 

\section{Final comments and open questions}
\label{sec:discussion}

We first show that several bijections are underlying the results
presented in this paper. We then raise a number of questions of an
enumerative or graph-theoretic nature.

\subsection{Bijections}
In what follows, we discuss three objects closely related to the graph
$G_\pi$: first the graph itself,  second its oriented version
$\vec G_\pi$
(each edge is oriented from the vertex with the lower label to the
vertex with the higher label), and
finally its natural embedding in $\N^2$ (where the vertex $i$ is placed at
position $(i, \pi(i))$ and the edges are represented by
straight lines,  as on the left of
Figure~\ref{fig:definitions}). 

\subsubsection{The graph $G_\pi$}
We first note that the map $\pi \mapsto G_\pi$ is injective. That is,
one can recover $\pi$ from $G_\pi$. To see this, orient $G_\pi$ to obtain
$\vec G_\pi$. Then, for every
vertex $i$ in $\vec G_\pi$, 
let $a(i)$ be the number of vertices that can be 
reached from $i$ by a directed path. This is the number of $j\ge i$ such
that $\pi(j)\ge \pi(i)$, and the sequence $\pi(1), \pi(2), \ldots ,
\pi(n)$ can be easily reconstructed (in this order) from 
the list $(a(1), \ldots, a(n))$. For instance, if $m$ vertices can be reached
from $1$ (that is, $a(1)=m$), then it means that $ \pi(1)=n-m+1$, and so on (we have
assumed implicitly that $|\pi|=n$).

\begin{figure}[h]
\begin{center}
  \input{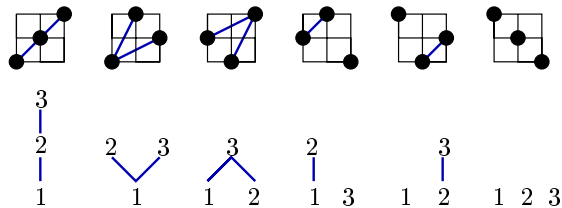}
\end{center}
\caption{The 6 posets obtained from \ps\ of length 3.}
\label{fig:posets}
\end{figure}

As noted at the beginning of the paper, $\vec G_\pi$ is the Hasse
diagram of a certain poset $P$ on $\llbracket n\rrbracket=\{1, 2, \ldots, n\}$. The underlying
order is \emm natural,, meaning that if $i<j$ in $P$, then $i<j$ in
$\N$. (We refer to~\cite[Chap.~3]{stanley-vol1} for generalities on
posets.) The $n!$ \ps\ of $\Sn_n$ thus provide $n!$ distinct natural
orders on $\llbracket n\rrbracket$. Not all natural orders are
obtained in that way: even for $n=3$, there are 7 natural orders but
only 6 \ps, and the poset in which the only relation is $1<3$ is not
obtained from any \p\ (Figure~\ref{fig:posets}). The posets that \emm are,
actually obtained from the construction $\pi 
\mapsto \vec G_\pi$ are, by definition, the natural orders on $\llbracket n\rrbracket$ \emm of dimension,
2~\cite[Exercise~3.10]{stanley-vol1}.
 
Some graph properties of $G_\pi$ easily follow from the
construction. For instance,
$G_\pi$ is isomorphic to $G_{\pi^{-1}}$ (more precisely,
$G_{\pi^{-1}}$ is obtained by relabelling the vertex $i$ by $\pi(i)$).
The natural embedding of  $G_{\pi^{-1}}$ is obtained by reflecting the
embedding of $G_\pi$ through the main diagonal. 
Of course, $G_\pi$ is triangle free (every Hasse diagram is). In
particular, $e(\pi)\leq\lfloor n^2/4\rfloor$ (see~\cite{aigner}) and 
it is easy to construct a permutation showing that this bound is tight. 
The number of edges of $G_\pi$ can also be interpreted in terms of
Bruhat order (see \cite[Exercise~3.75]{stanley-vol1}, \cite[Ch.~2]{bjorner}):  
 it is the  number of permutations covering (poset-wise) $\pi$ in the Bruhat
 order.  

\subsubsection{Rooted \tl\ \ps}
\label{sec:plane}
\begin{figure}[htb]
\begin{center}
  \input{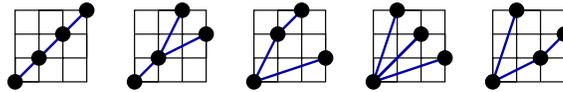}
\end{center}
\caption{The 5 rooted \tl \ \ps \ of length 4 and the corresponding plane trees.}
\label{fig:plane-trees}
\end{figure}

Here, we want to show a simple bijection between rooted
\tl\ \ps\ of size $n$ and plane trees with $n-1$ edges.
This explains why such \ps\ are counted by the  Catalan
number $C_{n-1}$. Recall
that \ps\ $\pi$ avoiding $21\bar 354$ are exactly those such that the
natural embedding of $G_\pi$ is planar (see the remark following
Lemma~\ref{lem:pattern}). This holds in particular for 
rooted \tl\ \ps: the embedding of $G_\pi$ is thus a (rooted) plane
tree.  Then, observe that the decomposition of \fl\ \ps\
illustrated in Figure~\ref{fig:structure}, once specialized to rooted \ps,
coincides with the standard decomposition of plane trees (a left
subtree joined to the root by an edge, and another plane tree, see
Figure~\ref{fig:plane}). This 
means that every plane tree is obtained from exactly one  rooted \tl
\ \p. This is illustrated in Figure~\ref{fig:plane-trees} for  \ps\ 
of length 4. 
 
\begin{figure}[t]
\begin{center}
  \input{plane-decomp1.pstex_t}
\end{center}
\caption{The decomposition of rooted \tl \ \ps.}
\label{fig:plane}
\end{figure}

\subsubsection{Path-like \ps}
\begin{figure}[htb]
\begin{center}
  \input{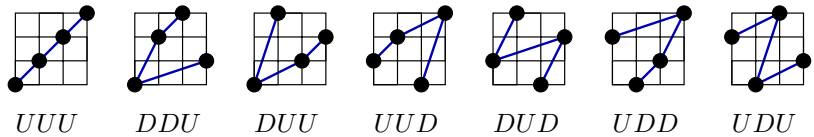}
\end{center}
\caption{The 7 path-like \ \ps \ of length 4 and the corresponding words.}
\label{fig:path}
\end{figure}

Consider a path-like \p\ $\pi$ of length at least 2. The graph $G_\pi$ has two vertices of
degree 1. Define a word $W(\pi)$  on the alphabet $\{U,D\}$ by
following the path $G_\pi$ from the vertex of degree 1 with the lowest label
to the other vertex of degree $1$, encoding each edge 
of this path by a letter $U$ (like
\emm up,) or $D$ (like \emm down,) depending on how the labels of the
vertices vary along this edge. Examples are shown in Figure~\ref{fig:path}.
It turns out that the map $W$ is a bijection from path-like  
\ps\ of length $n$ to words of length $n-1$ \emm distinct from,
$D^{n-1}$. In particular, this explains why the number of path-like
\ps\ of length $n$ is $2^{n-1}-1$. 

Again, this result follows from the decomposition of
path-like \ps\ that led to the equation of
Proposition~\ref{prop:bivariate}. Indeed, this decomposition gives, for the
\emm noncommutative, \gf\ defined by
$$
\mathcal{P}= \sum_{\pi \mbox{\tiny{\ path-like}}} W(\pi)
$$
the following equation:
$$
\mathcal{P}= \epsilon + U^+ + D^+U^+ + (\mathcal{P}-\epsilon)DU^*,
$$
where $\epsilon$ denotes the empty word and we have used the standard
notation $D^+=\sum_{i\ge 1}D^i$ and 
$U^*=\sum_{i\ge 0} U^i$. It is easy to see that the solution of this
equation is 
$$
 \mathcal{P}= \{U,D\} ^* - D^+.
$$
That is to say, the non-empty words $W(\pi)$ are those containing at
least one $U$, and each such word corresponds to a unique path-like
\p.

\bigskip

\subsection{Open problems}

\subsubsection{Enumeration}
In this paper, we have characterized and counted \fl\ \ps\ and some
of their natural subclasses. This work raises similar questions for
several supersets of \fl\ \ps. The most natural ones are probably the
following two:
\begin{enumerate}
  \item what is the number of \emm plane, \ps\ of $\Sn_n$, that is to
  say, \ps\ avoiding $21\bar354$?
\item what is the number of \ps\ associated with a \emm Gorenstein,
  Shubert variety? These \ps\ generalize \fl\ \ps, and have been
  characterized in~\cite{WooYong1}. 
\end{enumerate}
We also recall that the enumeration of 1324 avoiding \ps\ is still an
open problem~\cite{rechni,Marinov}. Permutations avoiding 2143 are
called \emm vexillary, and are equinumerous with 1234 avoiding
\ps\ \cite{bwx,west-these}, 
which have been enumerated in~\cite{gessel-symmetric}.

\smallskip
Another natural question is to count \ps \ $\pi$ by their length and
the number $e(\pi)$  of bars in their bar diagram (which is the number
of \ps\ covering $\pi$ in the Bruhat order). To our knowledge, the
bivariate series
$$
E(t,x)=\sum_{n\ge 0} \frac {t^n}{n!} \sum_{\pi \in \Sn_n} x^{e(\pi)}
$$
is not known. However, the \emm total, number of edges in the bar
diagrams of \ps\ of $\Sn_n$ \emm is, known: if 
$$
e(n) =\sum_{\pi \in \Sn_n} {e(\pi)},
$$
then $$e(n)
= (n+1)! (H(n+1)-2) +n!
$$
where $H(n)=1+1/2+ \cdots + 1/n$ is the $n$th harmonic number. Indeed,
as communicated to us by David Callan, it is not hard to see that the number
of \ps\ of $\Sn_n$ having a bar going from $i$ to $j$, with $i<j$, is
$n!/(j-i+1)$, and the above result follows easily.

Note
that $e(n)$ is also the number of edges in the Hasse diagram of the
Bruhat order of $\Sn_n$. The exponential \gf\ of the numbers $e(n)$ is
$$
\sum_{n\ge 0} e(n) \frac {t^n}{n!} = \frac{\partial E}{\partial x}
(t,1)= \frac 1{(1-t)^2} \left( \log \frac 1{1-t} -t\right).
$$
The average number of bars in a \p\ of $\Sn_n$ is
$$ 
\frac{e(n)}{n!} =\log  \left( n \right) n+ \left( -2+\gamma \right) n+\log  \left( n
 \right) +1/2+\gamma+O ( 1/n)
$$
where $\gamma$ is Euler's constant. This can be compared to the
average number of non-inversions, which 
is known to be $n(n+1)/4$. 
%
Related questions have recently been studied in~\cite{adin}.

\subsubsection{Graph questions}
We have seen that the labeled  graphs obtained from the map
$\pi\mapsto G_\pi$ are the Hasse diagrams of natural orders of
dimension 2. 
One can also wonder which \emm unlabelled, graphs are obtained
through our construction. Clearly, these graphs must be triangle free.
  However, this is
not a sufficient condition. For example, 
by an exhaustive computer search one can verify
that the triangle-free graph formed
of the vertices and edges of a cube 
is not produced from any  permutation in $\Sn_8$.

Note that, by Section~\ref{sec:plane}, all unlabelled trees (and thus
all unlabelled  forests) are obtained through our construction.

\subsubsection*{Acknowledgements}
The authors would like to thank David Callan, Aur\'elie Cortez, Alexander Woo and Alexander Yong for their advice
%
%
and also the anonymous referees whose comments on a former draft of
this paper greatly improved its quality. 

\bibliographystyle{plain}
\bibliography{biblio.bib}

\end{document}